\documentclass[]{gOMS2e}
\usepackage{float,latexsym,graphicx,shortvrb}
\usepackage{amsmath,amsfonts,bm,mathtools,booktabs}
\usepackage{amsmath,amsfonts,bm,mathtools,graphicx,latexsym,shortvrb}
\usepackage{xcolor}
\usepackage{hyperref}
\usepackage{multirow}
\usepackage{colortbl}

\def\btheo{\begin{theorem}}
\def\etheo{\end{theorem}}
\def\bprop{\begin{proposition}}
\def\eprop{\end{proposition}}
\def\bexam{\begin{example}}
\def\eexam{\end{example}}
\def\bdefi{\begin{definition}}
\def\edefi{\end{definition}}
\def\blemm{\begin{lemma}}
\def\elemm{\end{lemma}}

\newcommand{\st}{\textrm{s.t.}}
\newcommand{\sta}{\textrm{sta}}
\def\inv{{-1}}
\def\Col{\textrm{\em C$_{ol}$}}

\def\det{\textrm{det }}

\def\non{\nonumber}

\def\[#1\]{\begin{align}#1\end{align}}
\def\bcase{\begin{cases}}
\def\ecase{\end{cases}}
\def\bpmat{\begin{pmatrix}}
\def\epmat{\end{pmatrix}}
\def\bbmat{\begin{bmatrix}}
\def\ebmat{\end{bmatrix}}
\def\beqn{\begin{eqnarray}}
\def\eeqn{\end{eqnarray}}
\def\beqnx{\begin{eqnarray*}}
\def\eeqnx{\end{eqnarray*}}
\def\beq{\begin{equation}}
\def\eeq{\end{equation}}
\def\bitem{\begin{itemize}}
\def\eitem{\end{itemize}}
\def\btheo{\begin{theorem}}
\def\etheo{\end{theorem}}
\def\bblock{\begin{block}}
\def\eblock{\end{block}}
\def\benum{\begin{enumerate}}
\def\eenum{\end{enumerate}}

\def\bxx{\bar{\bm x}}

\def\txx{\tilde{\bm x}}
\def\bary{\bar{y}}
\def\byy{\bar{\bm y}}

\def\hatf{\hat{f}}
\def\hff{\hat{\bm f}}
\def\hatp{\hat{p}}

\def\hpp{\hat{\bm p}}

\def\xx{\bm x}
\def\yy{\bm y}
\def\zz{\bm z}

\def\ff{\bm f}

\def\pp{\bm p}

\def\tt{\bm t}

\def\II{\bm I}

\def\GG{\bm G}
\def\QQ{\bm Q}

\def\UU{\bm U}

\def\a{\alpha}

\def\ve{\varepsilon}

\def\l{\lambda}

\def\s{\sigma}

\def\t{\tau}

\def\L{\Lambda}

\def\bgs{\bar{\s}}

\def\tldgl{\tilde{\l}}
\def\bgmu{\bar{\mu}}

\def\bbR{\mathbb{R}}

\def\calP{\mathcal{P}}

\def\calS{\mathcal{S}}

\def\calX{\mathcal{X}}

\def\lpa{\left(}
\def\rpa{\right)}

\def\lbc{\left\{}
\def\rbc{\right\}}

\def\half{\frac{1}{2}}

\newcommand{\qed}{\hspace*{\fill} $\Box$ \vspace{2ex}}
\newcommand{\real}{\mathbb{R}}

\begin{document}
\doi{10.1080/1055.6788.YYYY.xxxxxx}
\issn{1029-4937}
\issnp{1055-6788}
\jvol{00} \jnum{00} \jyear{2013} \jmonth{June}

\markboth{Taylor \& Francis and I.T. Consultant}{Optimization Methods and Software}

\articletype{}

\title{\bf Global Solutions to Large-Scale Spherical Constrained  Quadratic Minimization via Canonical Dual Approach}

\author{Yi Chen$^\dagger$ and David Y. Gao$^\dagger$$^{\ast}$\thanks{$^\ast$Corresponding author. Email: d.gao@ballarat.edu.au
\vspace{6pt}}
\\\vspace{6pt}  $^\dagger${\em{School of Science, Information Technology and Engineering, University of Ballarat, Victoria, 3353, Australia}}\\\vspace{6pt}\received{} }

\maketitle

\begin{abstract}
This paper presents
 global optimal solutions to a  nonconvex quadratic minimization problem  over a sphere constraint.
 The problem is well-known as a trust region subproblem and has been studied extensively
  for decades. The main challenge is   the so called 'hard case', i.e., the problem has multiple solutions on the  boundary  of the sphere.
  By canonical duality theory, this challenging problem is able to reformed as an one-dimensional canonical
  dual problem  without duality gap.
    Sufficient and necessary conditions are obtained by the triality theory, which can be used to identify
     whether the problem is hard case or not.
     A perturbation method and the associated algorithms are proposed to solve this hard case problem.
     Theoretical results and methods are verified by large-size examples.   \bigskip

\begin{keywords} global  optimization; quadratic minimization problems;
 canonical duality theory; trust region subproblem
\end{keywords}
\begin{classcode}90C20; 	90C26; 90C46
\end{classcode}\bigskip

\end{abstract}

\section{Introduction}\label{se:intro}

We consider the following  quadratic minimization problem:
\[
(\calP)~~\min &~~ P(\xx)=\xx^T\QQ\xx-2\ff^T\xx \label{p:prim}\\
\st &~~ \xx\in\calX , \non
\]
where  the given matrix $\QQ\in\bbR^{n\times n}$ is assumed to be  symmetric,  $\ff \in \bbR^n$ is an arbitrarily given vector,
 and the feasible region is defined as
\[
\calX=\lbc \xx\in\bbR^n~|~\|\xx\|\leq r\rbc,\label{eq:domain}
\]
in which,  $r$ is a positive real number.

 Problem  $(\calP)$ arises naturally in computational mathematical physics with extensive applications in  engineering sciences.
From the point view of systems theory, if the vector $\ff \in \real^n$ is considered as
an  input (or source), then the solution $\xx \in \real^n$ is refereed as the output (or state) of the system.
By the fact that the capacity of  any given  system is  limited,   the spherical constraint in
$\calX$ is naturally required for virtually every real-world system.
For example, in engineering structural analysis, if the applied force field $\ff \in \real^\infty$
is  big enough, the stress distribution
 in the structure will reach its elastic limit and the structure will  collapse.
For  elasto-perfectly plastic materials, the  well-known
 von Mises yield condition is  a quadratic inequality constraint at each material point\footnote{The Tresca yield condition
 is equivalent to a box constraint at each material point}
 (see Chapter 7, \cite{Gao00duality}). By finite element method,  the variational problem  in structural limit analysis
can be formulated as a large-size nonlinear optimization  problem with $m$ quadratic inequality constraints (the $m$ depends on the
number of total finite elements). Such problems have been studied extensively in
computational mechanics for more than fifty years and the so-called penalty-duality finite element programming \cite{gao-cs-88,gao-ijss-88}
is one of well-developed efficient methods for
solving this type of problems in engineering sciences.

In mathematical programming, the problem $(\calP)$
is  known  as a trust region subproblem, which arises in trust region methods \cite{powell03trust}.
A more general problem with nonconvex  quadratic constraint  is considered in \cite{Xing13qpqc}.
Although the function $P(\xx)$ can be nonconvex if the matrix $\QQ$ has negative eigenvalues,  it is proved that the problem ($\calP$) is
hidden convex, i.e.  ($\calP$) is actually equivalent to a convex optimization problem \cite{BenTal96hidden}.
By the optimization theory we know
 that the vector $\bxx$ is a solution of  ($\calP$)  if there exists a Lagrange multiplier   $\bgmu$ such that the following conditions hold \cite{Flippo96QOball}:
\[
&(\QQ+\bgmu \II)\bxx=\ff\label{eq:kkt1} \\
&\|\bxx\|\leq r\label{eq:kkt2}\\
&\QQ+\bgmu \II\succeq0,~\bgmu\geq0\label{eq:kkt3}\\
&\bgmu(\|\bxx\|-r)=0\label{eq:kkt4}
\]

Let $\l_1$ be the smallest eigenvalue of the matrix $\QQ$. From conditions (\ref{eq:kkt3}), we know that
$$
\bgmu\geq\max\{0,-\l_1\}.
$$
If the problem ($\calP$) has no solution on the boundary of $\calX$,
 then $\QQ$ must be positive definite and $\|\QQ^\inv\ff\|<r$, which leads to $\bgmu=0$.
  If  ($\calP$) has a solution on the boundary of $\calX$ and $(\QQ+\bgmu \II)\succ 0$,
  then we have  $\|(\QQ+\bgmu \II)^\inv\ff\|=r$. In this case,  the multiplier $\bgmu$ can be found by using Newton's method.
 However, if the  solution $\bxx$ is located on the boundary of $\calX$ and  $\det (\QQ+\bgmu \II) = 0$,
  this situation is the so-called `hard case' \cite{more83QOball},
    which leads to numerical difficulties.
 In this case, the equation $(\QQ+\bgmu \II)\xx=\ff$ has no unique solution,
 and all vectors in the form $\xx=(\QQ+\bgmu \II)^\dagger\ff+\t\txx$ with $(\QQ+\bgmu \II)\txx=0$ are its solutions.
  In \cite{more83QOball}, Mor\'e and Sorensen proposed a safeguarding scheme to update $\mu$ and replaced $\txx$ by the vector $\zz$ with $\|R\zz\|$ being an approximation of the smallest singular value of $R$, where $R$ is the Cholesky factorisation of $\QQ+\mu \II$.
    Many other methods have been developed to deal with either  hard case or  large-size problems. Methods through a
    parameterized eigenvalue problem are discussed in \cite{sorensen97QOball,Rendl97QPballSDP,rojas01trs,fortin04QOball}.
     At each iteration, the Lanczos method was used to calculate an approximation of the smallest eigenvalue. Another kind of methods \cite{steihaug83trs,gould99QOball,Hager01QOball} searches solutions in the Krylov space, which is gradually expanded during iterations. In \cite{tao98DCtrs}, the d.c. (difference of convex functions) algorithm is applied to solve the problem ($\calP$).

The  goal of this paper is to solve the  problem ($\calP$) in any size, especially for the hard case.
Our approach is the canonical duality theory,   a newly developed and potentially powerful methodological theory,
which  has been used successfully for solving a
large class of nonconvex/nonsmooth/discrete problems in analysis and global optimization within a unified framework (see
\cite{Gao09unified,Gao09Lagrangian}).
 This theory  is composed mainly of
(1) a canonical dual transformation; (2) a complementary-dual principle, and (3) a triality theory.
We first show in the next section  that by the canonical dual transformation, this constrained nonconvex problem can be reformed as
  a one-dimensional optimization problem. The complementary-dual principle shows that this one-dimensional problem
   is canonically (i.e. perfectly) dual to $(\calP)$ in the sense that both problems have the same set of KKT solutions.
While the triality theory (mainly the first statement, i.e. the canonical min-max duality) provides
 sufficient and necessary conditions  for identifying  global optimal solutions.
 In order to solve the hard case,   a perturbation method is proposed in  Section \ref{se:pert} and, accordingly,
 a canonical primal-dual
  algorithm is developed in Section \ref{se:algorithm}. Numerical results  presented in Section \ref{se:nume} show that our approach can
  efficiently solve  large-size   problems.
 The paper is ended with some conclusion remarks.

\section{Canonical dual problem}\label{se:cano}
According to  \cite{Gao09Lagrangian}, the canonical dual problem of $(\calP)$ is given by
\[
\sta\lbc P^d(\s) ~|~ \s\in\calS_a\rbc , \label{p:dual}
\]
where the notation $\sta$ denotes computing stationary points of the canonical dual function  $P^d(\s)$  which is defined as
\[
P^d(\s)=-\ff^T\GG_a(\s)^{-1}\ff-r^2\s,\label{ep:dualfun}
\]
in which,  $\GG_a(\s)=\QQ+\s \II$ and $\GG_a(\s)^{-1}$ denotes the inverse of $\GG_a(\s)$. The  feasible set $\calS_a$ is   defined as
$$\calS_a=\{\s ~|~\s\geq0,~ \ff\in\Col(\GG_a(\s))\},$$
and the notation  $\Col(\cdot)$ represents the   column space of $\GG_a$.

We note that the canonical dual   $P^d(\s)$ is a function of a scalar variable $\s \in \real$, regardless of
 the dimension of the primal problem.
The canonical duality theory demonstrates that there is no duality gap between the primal problem ($\calP$) and its  canonical dual
 (\ref{p:dual}), which is illustrated by the following theorem.

\begin{theorem}\label{th:AnalSolu}
{\rm (\textbf{Analytical Solution and Complementary-Dual Principle} \cite{Gao00duality,Gao09Lagrangian})}
The problem (\ref{p:dual}) is canonically dual to the problem ($\calP$) in the sense that if $\bar{\s}\in\calS_a$ is a critical point of $P^d(\s)$, then
\beq\label{eq:solvedx}
\bar{\xx}=\GG_a(\bar{\s})^{-1}\ff
\eeq
is a KKT point of the primal problem ($\calP$), and we have
\beq\label{eq:nogap}
P(\bar{\xx})=P^d(\bar{\s}).
\eeq
\end{theorem}

 The proof is omitted here, which is analogous with that in \cite{Gao09Lagrangian}.
In order to identify global optimal solutions among all the critical points of $P^d(\s)$, a subset of $\calS_a$  is needed:
\[
\calS_a^+=\lbc \s\in\calS_a ~|~ \GG_a(\s)\succeq\bf{0} \rbc.\non
\]
Therefore, the canonical dual problem of ($\calP$) can be proposed as the following
\[\label{p:dualSplus}
(\calP^d)~~\max \{P^d(\s)~|~ \s\in\calS_a^+\}.
\]

\begin{theorem}\label{th:global}
{\rm (\textbf{Global Optimality Condition} \cite{Gao00duality,Gao09Lagrangian})}
Suppose that $\bar{\s}$ is a critical point of $P^d(\s)$. If $\bar{\s} \in\calS_a^+$ and $\det(\GG_a(\bar{\s}))\neq 0$, then $\bar{\s}$ is
 a global maximal solution of the problem ($\calP^d$) on $\calS^+_a$
 and $\bar{\xx} =\GG_a(\bar{\s})^{-1}\ff  $   is a  global minimal solution of the primal problem ($\calP$), i.e.
\[
P( \bar{\xx} )=\min_{\xx\in\calX} P(\xx)=\max_{\s\in\calS_a^+}P^d(\s)=P^d( \bar{\s}).\label{eq:eqglobal}
\]
\end{theorem}

According to  the triality theorem \cite{Gao00duality,gao12triality},
the global optimality condition (\ref{eq:eqglobal})  is called canonical min-max duality.
 It guarantees that if there is a critical point in the interior of $\calS_a^+$,
  computing the global minimal solution of the nonconvex problem ($\calP$)
 can be converted to a   concave maximization problem.
 Therefore, the so-called hidden convexity discovered in \cite{BenTal96hidden} is actually a special case of the
 canonical duality. Also, by the canonical duality theory,   complete solutions of the problem ($\calP$) have been discussed by Gao in
  \cite{Gao04quadratic}, wherein, Theorem 3  states that if $P(\xx)$ is not  convex and the Morse index of $P(\xx)$
  (i.e. the number of negative eigenvalues
   of $\QQ$, see Chapter 5 in \cite{Gao00duality}) is $i_d$,
   then the problem ($\calP$) has at most $2i_d+1$ KKT points on the boundary of the sphere and they can be calculated from the $2i_d+1$ KKT points of the dual problem. Moreover, the corresponding primal and dual functions    are equal at each of these KKT points.
    This theorem presents a perfect duality relationship between the problem ($\calP$) and its dual problem.
    By the
    double-min duality statement in the weak-triality theory proven recently
    (see \cite{gao12triality,morelas-gao-naco,Morales11triality}), we know that the
    problem ($\calP$)  has at most one local  minimizer since
    the  canonical dual problem is in one-dimensional space.
    Similar result is also proven in \cite{martinez94qpball}.
For the hard case, i.e. the matrix $ \GG_a(  {\s} )  $ is singular at the critical point  $\bar{\s}$,
the canonical dual $P^d(\s)$ should be replaced by (see \cite{Gao10boxIP})
\[
P^d(\s)=-\ff^T \GG_a(\s)^\dag \ff -r^2\s , \label{ep:dualfung}
\]
where $\GG_a(\s)^\dag $ stands for a generalized inverse of $\GG_a(\s)$.
Since this function is not strictly concave on  $\calS_a^+$,
it may have multiple critical points located on the boundary of $\calS_a^+$.
In the following sections, we will first study the existence conditions of these critical points, and then to study associated algorithm  for computing these solutions.

\section{Existence condition}\label{se:exist}

By the symmetry of  the matrix $\QQ$,
there exist diagonal matrix $\L$ and orthogonal matrix $\UU$ such that $\QQ=\UU\L \UU^T$.
 The diagonal entities of $\L$ are the eigenvalues of the matrix $\QQ$ and are arranged in nondecreasing order,
$$
\l_1=\cdots=\l_k<\l_{k+1}\leq\cdots\leq\l_n.
$$
The columns of $\UU$ are corresponding eigenvectors.

%

It's easy to verify that $(\QQ+\s \II)^{-1}=\UU(\L+\s \II)^{-1}\UU^T$. Let $\hff=\UU^T\ff$. Therefore, we can rewrite the dual function into
\[
P^d(\s)&=-\frac{\sum_{i=1}^k\hatf_i^2}{\l_1+\s}-\sum_{i=k+1}^n\frac{\hatf_i^2}{\l_i+\s}-r^2\s,\label{eq:diagPd}
\]
where $\hatf_i$ are elements of $\hff$. We notice that the function $P^d(\s)$ is always well defined and have stationary points over its domain except that $\ff=0$.
Thus, for the case of  $\ff\neq0$, the dual problem ($\calP^d$) is well defined. For the case of  $\ff=0$, the canonical dual problem can be solved by the perturbation method provided in the next section.

\begin{proposition}\label{th:exist}
{\rm (\textbf{Existence Condition})}
Suppose that $\l_i$ and $\hatf_i$ are defined as above and there is a solution of the problem ($\calP$) on the boundary of $\calX$.
Then there exists a critical point of $P^d(\s)$ in $(-\l_1,+\infty)$ if and only if either
$\sum_{i=1}^k\hatf_i^2\neq0$ or $\sum_{i=k+1}^n\frac{\hatf_i^2}{(\l_i-\l_1)^2}>r^2$.
 If $P^d(\s)$ has a critical point  $\bar{\s}$ in $(-\l_1,+\infty)$, then this  critical point is unique
 and  $\bar{\xx}=\GG_a(\bar{\s})^\inv\ff$ is a global solution of the problem ($\calP$).
\end{proposition}

\proof:
First, let us  prove that  $P^d(\s)$ has  a critical point  in  $(-\l_1,+\infty)$ implies either  $\sum_{i=1}^k\hatf_i^2\neq0$ or $\sum_{i=k+1}^n\frac{\hatf_i^2}{(\l_i-\l_1)^2}>r^2$.
 Equivalently, we can prove that if $\sum_{i=1}^k\hatf_i^2=0$ and $\sum_{i=k+1}^n\frac{\hatf_i^2}{(\l_i-\l_1)^2}\leq r^2$
 the dual function $P^d(\s)$ will have no critical points in $(-\l_1,+\infty)$.
  If $\sum_{i=1}^k\hatf_i^2=0$, the first item vanishes in the expression (\ref{eq:diagPd}).
   Since we assume that $\sum_{i=k+1}^n\frac{\hatf_i^2}{(\l_i-\l_1)^2}\leq r^2$, the first-order derivative of the dual function $P^d(\s)$
\[
&(P^d(\s))'=\sum_{i=k+1}^n\frac{\hatf_i^2}{(\l_i+\s)^2}-r^2 \label{eq:der1}
\]
is always negative in  $(-\l_1,+\infty)$. Therefore, the dual function $P^d(\s)$  will have no critical points in $(-\l_1,+\infty)$.

Next we give the proof of the sufficiency, which is divided into two parts:

1) If $\sum_{i=1}^k\hatf_i^2\neq0$, $\l_1$ is a pole of $P^d(\s)$, which implies that as $\s$ approaches $-\l_1$
from the right side, the  function  $P^d(\s)$ approaches $-\infty$.
Also, $P^d(\s)$   approaches $-\infty$ as $\s$ approaches $+\infty$.
Therefore, $- P^d(\s)$ is coercive on $(-\l_1,+\infty)$.
Since, for any $\s\in(-\l_1,+\infty)$, $\GG_a(\s)$ is positive definite, $P^d(\s)$ is strictly concave on $(-\l_1,+\infty)$.
Thus there exists a unique critical point on $(-\l_1,+\infty)$.

2) If $\sum_{i=1}^k\hatf_i^2=0$ and $\sum_{i=k+1}^n\frac{\hatf_i^2}{(\l_i-\l_1)^2}>r^2$, $(P^d(\s))'$ is positive at $\s=-\l_1$.
Moreover, $(P^d(\s))'$ approaches $-r^2$ as $\s$ approaches $\infty$.
 Therefore, there exists at least one root for the equation $(P^d(\s))'=0$ over  $(-\l_1,+\infty)$,
  which means $P^d(\s)$ has at least one critical point in $(-\l_1,+\infty)$.
   Similarly, because of the strict concavity of $P^d(\s)$ over $(-\l_1,+\infty)$, the critical point is unique.

Let $\bar{\s}$ denote the critical point. If $\l_1\leq0$, we have $\bar{\s}\in\calS_a^+$. Then, from Theorem \ref{th:global},
 we further have that  $\bar{\xx}=\GG_a(\bar{\s})^\inv\ff$ is a global solution of the problem ($\calP$).
 If $\l_1>0$, from equations (\ref{eq:kkt1})--(\ref{eq:kkt4}), we know that the dual variable $\bgmu$
 satisfying $\|(\QQ+\bgmu \II)^\inv\ff\|=r$ is the critical point $\bar{\s}$.
 Thus $\bar{\xx}=\GG_a(\bar{\s})^\inv\ff=(\QQ+\bgmu \II)^\inv\ff$ is a global solution of the problem ($\calP$).

The proposition is proved.
\hfill\qed

\section{Perturbation method}\label{se:pert}
This section is devoted to the solutions for hard case, where
$\sum_{i=1}^k\hatf_i^2=0$ and $\sum_{i=k+1}^n\frac{\hatf_i^2}{(\l_i-\l_1)^2}\leq r^2$,
i.e. the existence condition obtained in the previous section
 is violated. This case leads to challenges for solving the problem  ($\calP$) via  (pure) mathematical analysis.
Our approach is the perturbation method, which has been used successfully in
canonical duality theory for solving nonlinear algebraic equations \cite{rgj},
 chaotic dynamical systems \cite{rgIMA}, as well as a class of NP-hard problems in global optimization \cite{Gao10boxIP,wang11maxcut}.
 In order to  reinforce the existence condition, a set of
 perturbation parameters
  \[
\a_i\neq0, ~~ \textrm{ for some } i\in\{1,\ldots,k\},\label{eq:constraints}
\]
   is introduced, and we let
$$
\pp=\ff+\sum_{i=1}^k\a_i\UU_i,~~
\hpp=\UU^T\pp, ~~
P_\a(\xx)=\xx^T\QQ\xx-2\pp^T\xx.
$$
Then the   perturbed problem  can be  defined as
\[
(\calP_\a)~~\min ~P_\a(\xx)~~~~\st~\xx\in\calX. \label{p:pert}
\]
It is true that the existence condition holds for the perturbed problem  since
  (\ref{eq:constraints}) will guarantee $\sum_{i=1}^k\hatp_i^2\neq0$.

The following theorem states that  for certain appropriate   $\{ \a_i\}_{i=1}^k$,  the optimal solution of the  perturbed problem
converges to that of the primal problem $(\calP)$.

\begin{theorem}\label{th:accuracy}
Suppose that $\l_1\leq 0$, and $\bxx$ and $\bxx^*$ are optimal solutions of the problems ($\calP$) and   ($\calP_\a$), respectively, on the boundary of $\calX$. Then, for any $\ve>0$, if the  parameters $\{\a_i\}$  satisfy
 \[
 \sum_{i=1}^k\a_i^2\leq (\l_2-\l_1)^2 \left(r^2-\sum_{i=k+1}^n\frac{\hatf_i^2}{(\l_i-\l_1)^2} \right) (1/\sqrt{2(1-\cos (\ve/r)}-1)^{-2} ,
 \]
 we have $\|\bxx^*-\bar{\xx}\|\leq\ve$.
\end{theorem}

\proof:
For simplicity, we rotate the coordinate system and  substitute $\xx$ with $\UU\yy$ in the problem ($\calP$).
As $\hatf_i=0$ for $i=1,\ldots,k$, variables $y_i$ for $i=1,\ldots,k$ appear in the form of squares in the target
 function. Since it is assumed that $\bxx$ and $\bxx^*$ are optimal solutions on the boundary of $\calX$,
both should satisfy the equality constraint in $\calX$.
Let  $\yy_{\ell} =\{y_i\}_{i=k+1}^n$.
On    the boundary of $\calX$,
the problem ($\calP$)
 is equivalent to the following  problem in $\real^{n-k}$:
\[\label{p:Pdiagl}
\min_{\|\yy_{\ell} \|\leq r} ~P^{\ell} (\yy_{\ell} )=\sum_{i=k+1}^n(\l_i-\l_1)y_i^2-\sum_{i=k+1}^n2\hatf_iy_i+\l_1r^2.
\]
Similarly, the  perturbed problem (\ref{p:pert}) with the equality constraint is equivalent to
\[\label{p:Pbdiag}
\min_{\|\yy_{\ell} \|\leq r} P_\a^{\ell} (\yy_{\ell} )=&\sum_{i=k+1}^n(\l_i-\l_1)y_i^2-\sum_{i=k+1}^n2\hatf_iy_i+\l_1r^2\non\\
&-2\sqrt{\sum_{i=1}^k\a_i^2}\sqrt{r^2-\sum_{i=k+1}^ny_i^2}.
\]
Then it is not difficult  to verify that $\byy_{\ell}=\{\bary_i=\UU_i^T\bxx\}_{i=k+1}^n$ and $\byy_{\ell} ^*=\{\bary_i^*=\UU_i^T\bxx^*\}_{i=k+1}^n$
are optimal solutions of problems (\ref{p:Pdiagl}) and (\ref{p:Pbdiag}), respectively.
Since $P^{\ell} (\yy_{\ell} )$ is a strictly convex function, it has a unique stationary point,
which is $\{\frac{\hatf_i}{\l_i-\l_1}\}_{i=k+1}^n$. Combining with the assumption,
we know that this stationary point is the global optimal solution of the problem (\ref{p:Pdiagl}), i.e.
$$
\bary_i=\frac{\hatf_i}{\l_i-\l_1}, ~i=k+1,\ldots,n.
$$
The function $P_\a^{\ell} (\yy_{\ell} )$ is also strictly convex. Furthermore, for any $\|\yy_{\ell} \|< r$, we have $P_\a^{\ell} (\yy_{\ell} )<P^{\ell} (\yy_{\ell} )$,
 and for any $\|\yy_{\ell} \|=r$, we have $P_\a^{\ell} (\yy_{\ell} )=P^{\ell} (\yy_{\ell} )$, which indicates that the unique stationary point of
 $P_\a^{\ell} (\yy_{\ell} )$ is in the interior of $\|\yy_{\ell} \|\leq r$. Thus it is the global optimal solution of the problem (\ref{p:Pbdiag}) and $\byy_{\ell} ^*$ satisfies
$$
\bary_i^*=\frac{\hatf_i}{\l_i-\l_1+\sqrt{\sum_{i=1}^k\a_i^2}(r^2-\sum_{i=k+1}^n(\bary_i^*)^2)^{-\half}}, i=k+1,\ldots,n.
$$
Obviously,
\[\label{eq:ybarlessy}
|\bary_i^*|<|\bary_i|, i=k+1,\ldots,n.
\]
We have the inequality
\[
\byy^{* T}\byy&=\sqrt{r^2-\byy_{\ell} ^{* T}\byy_{\ell} ^*}\sqrt{r^2-\byy_{\ell} ^{T}\byy_{\ell} }
+\byy^{*T}_{\ell} \byy_{\ell} \non\\
&\leq\half\lpa r^2-\byy_{\ell} ^{* T}\byy_{\ell} ^*+r^2-\byy_{\ell} ^{T}\byy_{\ell}  \rpa
+\byy^{* T}_{\ell} \byy_{\ell} \non\\
&=r^2-\half\|\byy_{\ell} ^*-\byy_{\ell} \|^2,\non
\]
which further implies that
\[\label{eq:disybary}
\|\byy^*-\byy\|\leq r\arccos \left(\frac{\byy^{* T}\byy}{r^2 } \right)
\leq r\arccos \left( \frac{r^2-\half\|\byy_{\ell} ^*-\byy_{\ell} \|^2}{r^2} \right) .
\]
Thus, if we want $\|\byy^*-\byy\|\leq \ve$, it is necessary to make sure $\|\byy_{\ell} ^*\byy_{\ell} \|^2\leq2r^2(1-\cos\frac{\ve}{r})$.
Because of the inequality
\[\label{eq:disybarlyl}
\|\byy_{\ell} ^*-\byy_{\ell}\|^2&\leq\frac{r^2}{\lpa(\l_2-\l_1)\sqrt{r^2-\byy_{\ell} ^{* T}\byy_{\ell} ^*}/\sqrt{\sum_{i=1}^k\a_i^2}+1\rpa^2},
\]
if we let its right side be less than or equal to $2r^2(1-\cos\frac{\ve}{r})$, we obtain
\[\label{eq:boundai2ystar}
\sum_{i=1}^k\a_i^2\leq\frac{(\l_2-\l_1)^2(r^2-\byy_{\ell} ^{*T}\byy_{\ell} ^*)}{(1/\sqrt{2(1-\cos\frac{\ve}{r})}-1)^2}.
\]
Hence, combining with relations in (\ref{eq:ybarlessy}), we can state that $\|\byy^*-\byy\|\leq \ve$ if the following inequality is true
\[\label{eq:boundai2}
\sum_{i=1}^k\a_i^2\leq\frac{(\l_2-\l_1)^2(r^2-\sum_{i=k+1}^n\frac{\hatf_i^2}{(\l_i-\l_1)^2})}{(1/\sqrt{2(1-\cos\frac{\ve}{r})}-1)^2}.
\]
Since $\|\bxx^*-\bxx\|=\|\byy^*-\byy\|$, the theorem is proved. \hfill\qed

Theorem \ref{th:accuracy} shows that with certain  proper parameters $\{ \a_i\}_{i=1}^k$,
 the existence condition is guaranteed for the perturbed problem such that the
 perturbation method can be used to solve the hard case.
 As we know that in hard case,  the primal problem $(\calP)$ may have multiple solutions  $\{ \bxx \}$
  on the boundary of the feasible region $\calX$.
 By the fact that the  perturbed problem ($\calP_\a$) is strictly convex in the neighborhood of $\bxx^*$
 and its global minimal  solution $\bxx^*$ will approach to
 one of these    $\{ \bxx \}$, depending on the parameters $\{ \a_i\}$.
From the projection theorem, we know that the nearest points to  ${\bxx}$ and $\bxx^*$  in the subspace spanned   by $\{\UU_1,\ldots,\UU_k\}$ are $\sum_{i=1}^k(\bxx^{T}\UU_i)\UU_i$ and $\sum_{i=1}^k(\bxx^{* T}\UU_i)\UU_i$, respectively,
which have the following relationship
\[\label{eq:projection}
\|\bxx^*-\sum_{i=1}^k(\bxx^{*T}\UU_i)\UU_i\|^2<\|{\bxx}-\sum_{i=1}^k(\bxx^{T}\UU_i)\UU_i\|^2.
\]
Therefore, the perturbed solution  $\bxx^*$ is closer to the subspace spanned  by $\{\UU_1,\ldots,\UU_k\}$ than the solution ${\bxx}$.

\section{Canonical primal-dual algorithm}\label{se:algorithm}
Based on the results in the previous section,
 we are ready to present an  algorithm.
 The Lanczos method is employed to compute approximately the smallest eigenvalue and the corresponding eigenvector, which  will be used
  to construct the safeguarding and perturbation.
  A canonical  primal-dual iterative scheme is introduced, which
   is matrix inverse free. The  essential cost of this algorithm  is only  the matrix-vector multiplication.

The key step of this algorithm is to solve the following perturbed  canonical dual   problem:
\[\label{eq:dualfunper}
(\calP^d_\a): \;\; \max \left\{ P_\a^d(\s)=-\pp^T\GG_a(\s)^{-1}\pp-r^2\s \;\; | \; \; \s \in \calS_a^+ \right\}
\]
Let $\psi(\s)$ be its first-order derivative, i.e.,
$$
\psi(\s)=(P_\a^d(\s))^\prime=\pp^T\GG_a(\s)^{-1}\GG_a(\s)^{-1}\pp-r^2.
$$
The critical point of $P_\a^d(\s)$ in $\calS_a^+$ is a  unique solution to
the equation $\psi(\s)=0$ in $\calS_a^+$.
Thus we need  to compute the zero of $\psi(\s)$ in $\calS_a^+$ to find the critical point. The first and second order derivatives of $\psi(\s)$ are
\[
&\psi^\prime(\s)=-2\pp^T\GG_a(\s)^{-1}\GG_a(\s)^{-1}\GG_a(\s)^{-1}\pp,\label{eq:psi1d}\\
&\psi^{\prime\prime}(\s)=6\pp^T\GG_a(\s)^{-1}\GG_a(\s)^{-1}\GG_a(\s)^{-1}\GG_a(\s)^{-1}\pp.\label{eq:psi2d}
\]
It is noticed that $\psi(\s)$ is strictly decreasing and strictly convex over $\calS_a^+$,
 $\psi(\s)$ will approach $-r^2$ as $\s$ approaches infinity and $\s=-\l_1$ is a pole of $\psi(\s)$.
Here we use the Lanczos method to compute an approximation of the smallest eigenvalue
 of $\QQ$ and  the corresponding eigenvector,
 denoted by $\tldgl_1$ and $\tilde{\UU}_1$, respectively. Clearly, we have $\|\tilde{\UU}_1\|=1$.

If $\tldgl_1>0$, we can conclude that $\l_1>0$ since $\tldgl_1$ is always smaller than $\l_1$. The dual feasible region is $\calS_a^+=[0,+\infty)$.
Thus, if $\psi(0)\leq0$, the maximiser of the dual function over $\calS_a^+$ is $\s=0$, and $\xx=\GG_a(0)^{-1}\ff$ is the global solution of
the primal problem ($\calP$). If $\psi(0)>0$, there exists a critical point in $\calS_a^+$, which is also the unique critical point in $(-\l_1,+\infty)$.

If $\tldgl_1\leq0$, we always intend to calculate the critical point in $(-\l_1,+\infty)$.
However, $(-\l_1,+\infty)$ may be not the $\calS_a^+$, because it is possible that $\l_1>0$, especially when a large error tolerance is chosen.
Thus, we should check whether the critical point is on the right side of 0. If not, the maximiser in $\calS_a^+$ should be $\s=0$.
When the problem is in the hard case, it is rather rational to choose $\a\tilde{\UU}_1$ with a proper scaling parameter $\a$ as a perturbation to $\ff$.

Although the perturbed canonical dual problem $(\calP^d_\a)$ is strictly concave on the closed domain
$\calS^+_a$, its derivative $\psi(\s)$ could be ill-conditioned when $\s$ approaches to the pole.
Therefore, instead of nonlinear optimization techniques,
a  bisection method is used to find the zero of $\psi(\s)$ in $(-\l_1,+\infty)$.
Each time, as a dual solution $\s>-\l_1$ is obtained, the value of  $\psi(\s)$ is calculated and checked to see if
 it is equal to zero. For moderate-size problems, it is possible to calculate $\GG_a(\s)^{-1}\pp$ by computing the inverse
  or decomposition of $\GG_a(\s)$, but it is not possible for very large-size problems,
 especially when the memory is very limited. One alternative approach is to solve the following strictly convex minimisation problem,
\[\label{p:xsgivens}
\min_{\xx\in\bbR^n}~\xx^T\GG_a(\s)\xx-2\pp^T\xx,
\]
whose optimal solution is $x=\GG_a(\s)^{-1}\pp$.
 Actually, during iterations, we do not need to calculate
  $\psi(\s)$ every time, especially for $\s$ being on the left side of the zero and close to the pole.
   We discover  that for a given  $\s$,
 the value of $\psi(\s)$ is equal to the optimal value of the following unconstrained concave maximization  problem
\[\label{p:zsgivens}
\max_{\zz\in\bbR^n}~-\zz^T\GG_a(\s)\GG_a(\s)\zz+2\pp^T\zz-r^2.
\]
By  the fact that the value of the target function will increase during the iterations,
 we can stop solving the problem (\ref{p:zsgivens}) if   the target  function is larger than
  a threshold and then claim that the $\s$ must be at the left side of the zero.
  Thus, the ill-condition  in  computing   $\psi(\s)$ as  $\s$ approaching  to the pole  can be prevented. 

An uncertainty interval should be initialized before the bisection method is applied, and it is used to
 safeguard that the interval contains the critical point. For the right end of the interval,
 any large enough number can be chosen. Actually, an upper bound of the critical point can be calculated,
 then it can be the right end of the uncertainty interval.
  Denote $\bgs^* \in (-\l_1,+\infty)$ be  the critical point of $P_\a^d(\s)$.
   From the definition of $\psi(\s)$, we get
$$
\frac{1}{(\l_1+\bgs^*)^2}\sum_{i=1}^n\hatp_i^2-r^2\geq0.
$$
Hence, $\sqrt{\sum_{i=1}^n\hatp_i^2}/r=\|\pp\|/r$ can be an upper bound.
However, the bound $\|\pp\|/r$ may be not tight. A practical way is to let $\s=-\l_1$ as a starting point and
  then to update $\s$   recursively by moving a certain step to its right each time.
   If the first $\s$ that makes the value of $\psi(\s)$ be negative is smaller than the upper bound $\|\pp\|/r$,
   it will be set to the right end of the uncertainty interval; otherwise, the upper bound will be the right end.

\begin{algorithm}[-- {\bf Initialization}]
\bitem
\item[]{\em Input: } coefficients $\QQ$, $\ff$ and $r$; a given error tolerance $\ve$.
\item[]{\em The smallest eigenvalue: } Use the Lanczos method to obtain $\tldgl_1$ and $\tilde{\UU}_1$.
\item[]{\em Perturbation: }
\bitem
\item[]{\bf If } the existence condition does not hold,   a perturbation is introduced and let
$$
\pp=\ff+\a\tilde{\UU}_1;
$$
\item[]{\bf Else } set $\pp=\ff$;
\item[]{\bf End if}
\eitem
\item[]{\em Uncertainty interval: } set an update size $s_t$ and a threshold $\ve_t$; let $\s=\s_{\ell} =-\tilde{\l}_1$;
\bitem
\item[]{\bf step 1:} Solve the problem (\ref{p:zsgivens}). If the value of the target  function is larger the  threshold $\ve_t$, the iteration stops, let $\s=\s+s_t$ and go to step 1; otherwise, go to step 2.
\item[]{\bf step 2:} Calculating the value of $\psi(\s)$.
\bitem
\item[] {\bf If } $\psi(\s)>0$, set $\s_{\ell} =\s$, $\s=\s+s_t$ and go to step 2;
\item[] {\bf Else } $\s_u=\s$ and STOP;
\item[] {\bf End if}
\eitem
\eitem
\eitem
\end{algorithm}

As the uncertainty interval $[\s_{\ell} ,\s_u]$ is obtained, the bisection method is applied to find the next iterate for $\s$, i.e. set $\s$ be the middle point of the  uncertainty interval. The main part of our algorithm is given as follows.
\begin{algorithm}[-- {\bf Main}]
\bitem
\item[]{\bf Do}
\bitem
\item[] set $\s=(\s_{\ell} +\s_u)/2$ and calculate the value of $\psi(\s)$;
\item[]{\bf If } $|\psi(\s)|<\ve$, then STOP and return $\s$ and $\xx$;
\item[] {\bf Else if } $\psi(\s)>0$, update $\s_{\ell} =\s$;
\item[] {\bf Else } update $\s_u=\s$;
\item[] {\bf End if}
\eitem
\item[]{\bf End do}
\eitem
\end{algorithm}

\section{Numerical experiments}\label{se:nume}
Let us  first present  two small-size examples
 to show the application of the canonical duality theory;
 we then list some large-size examples randomly generated  to demonstrate the efficiency of our method.

\subsection{Small-size examples}
\begin{figure}
\begin{center}
\subfigure[]{
\scalebox{0.45}[0.45]{\includegraphics{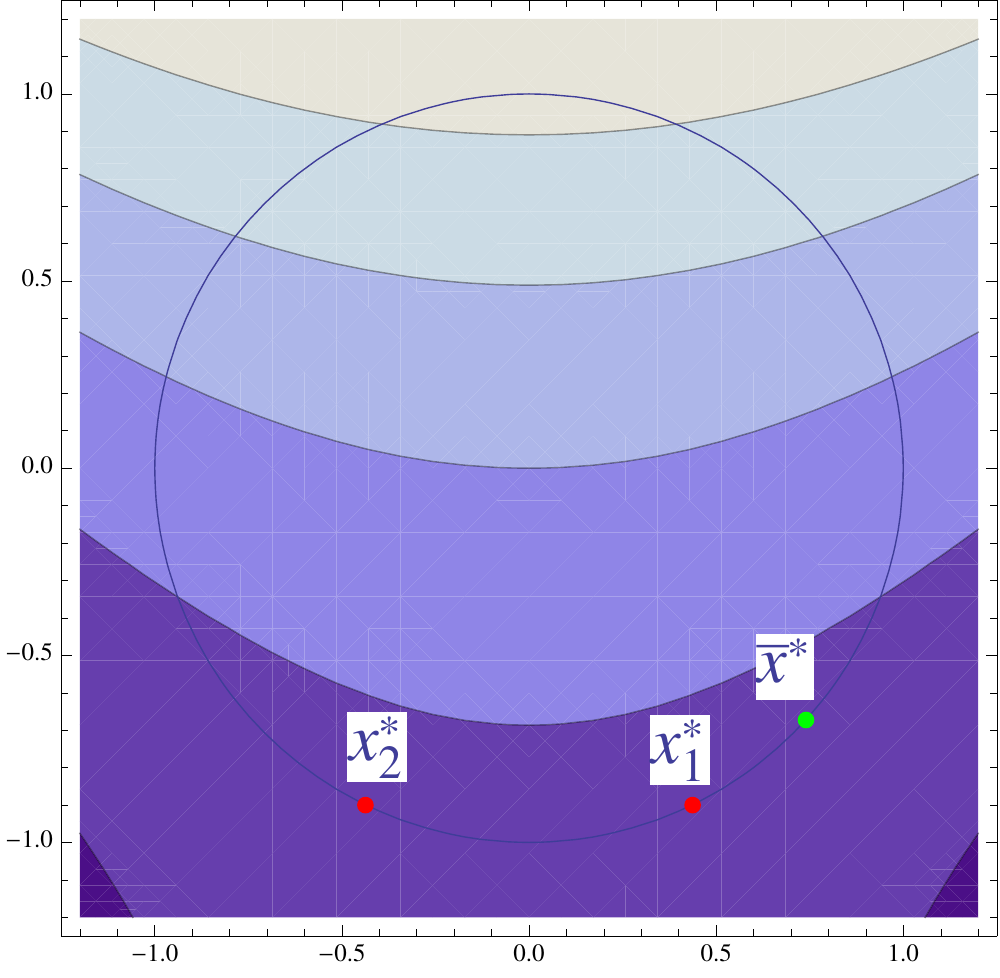}}}%
\subfigure[]{
\scalebox{0.65}[0.65]{\includegraphics{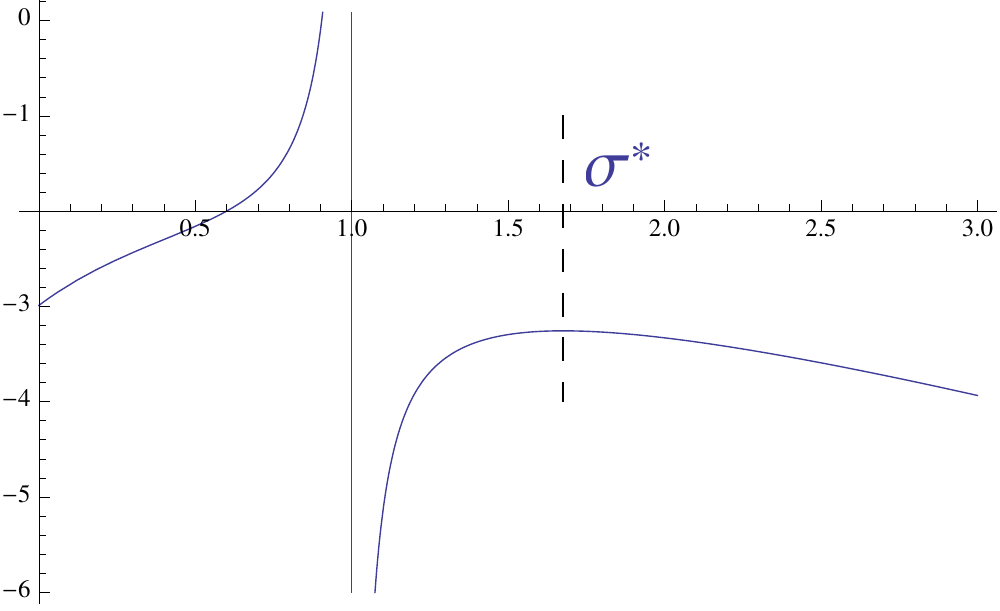}}}%
\caption{  $\pp=(0.5,-1.8)$. (a): The contours of the primal function and the boundary of the sphere; (b):   graph of the dual function.}%
\label{fig:ex1}
\end{center}
\end{figure}

\begin{figure}
\begin{center}
\subfigure[]{
\scalebox{0.45}[0.45]{\includegraphics{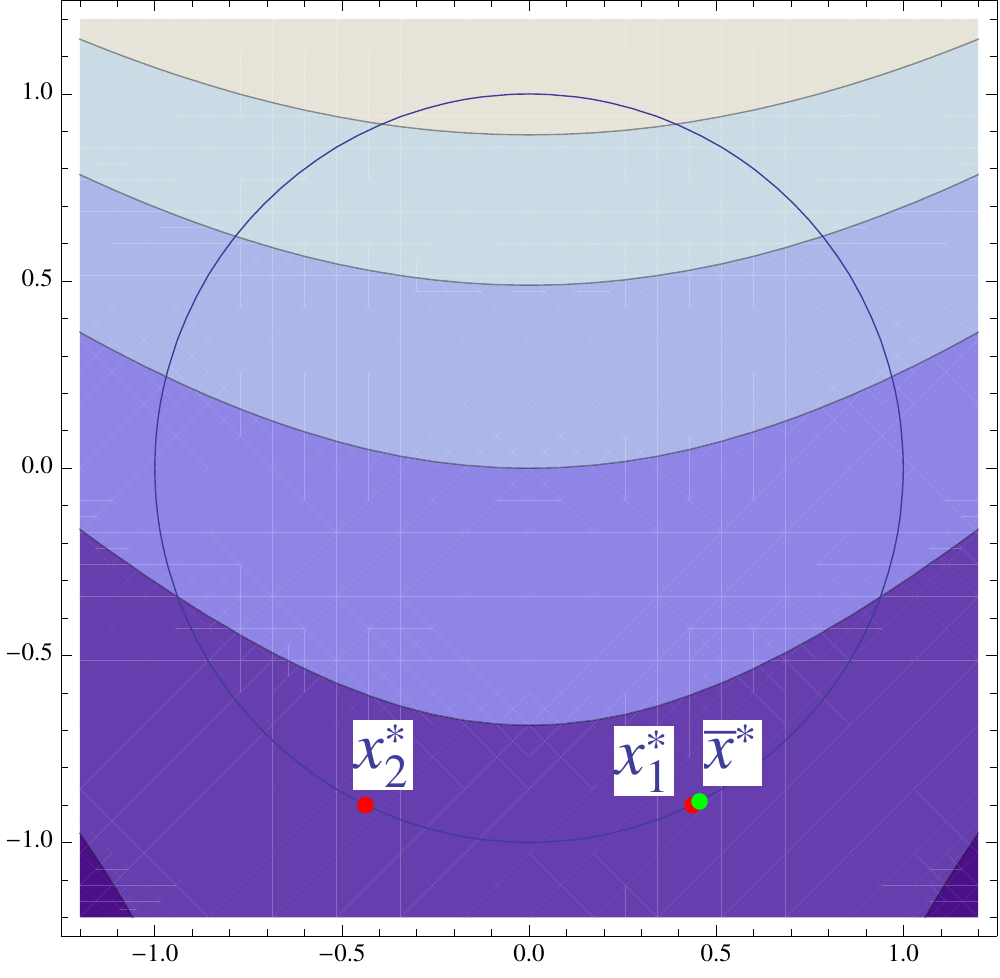}}}%
\subfigure[]{
\scalebox{0.65}[0.65]{\includegraphics{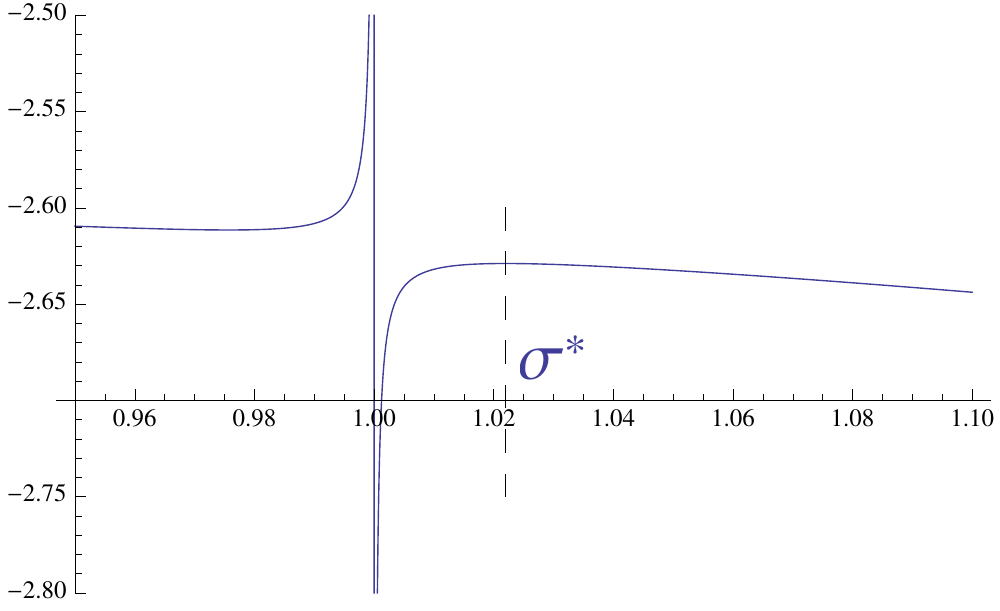}}}%
\caption{  $\pp=(0.01,-1.8)$. (a): The contours of the primal function and the boundary of the sphere; (b):   graph of the dual function.}%
\label{fig:ex1_2}
\end{center}
\end{figure}

{\bf Example 1} The given data  are
$$
\QQ=\bpmat -1 & 0\\0 & 1 \epmat,~~\ff=\bpmat 0\\-1.8 \epmat, ~~ r=1.
$$
The existence condition does not hold for this example. There are two global solutions, $\bxx_1=(0.437,-0.9)$ and $\bxx_2=(-0.437,-0.9)$, which are red points shown in Figure \ref{fig:ex1}.
In order to show how the perturbation method works, we first introduced a big perturbation in  the linear coefficient $\ff$ and
let $\pp=(0.5,-1.8)$.
The graph of the dual function of the perturbed problem is shown in the Figure \ref{fig:ex1}.
There is a critical point in the interior in $\calS_a^+$, which is $\bar{\s}=1.676$,
and the corresponding optimal solution for the perturbed problem is $\bxx^*=(0.74,-0.673)$, which is the green point in the Figure \ref{fig:ex1}.
  We then reduce the perturbation by letting $\pp=(0.01,-1.8)$.  The critical point is $\bar{\s}=1.022$ and the corresponding solution is $\bxx^*=(0.456,-0.89)$.
  Figure \ref{fig:ex1_2} shows that the perturbed solution $\bxx^*$ approaches $\bxx_1 $.

\begin{figure}
\begin{center}
\subfigure[]{
\scalebox{0.45}[0.45]{\includegraphics{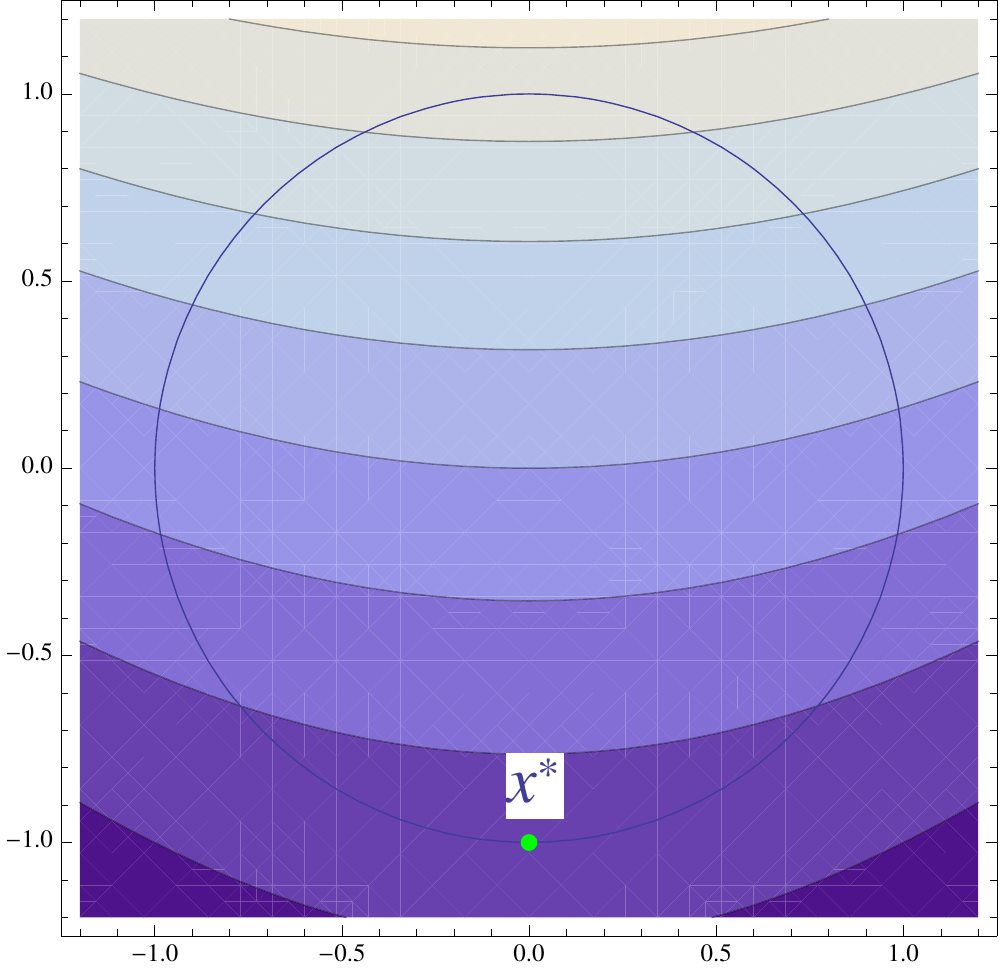}}\label{fig:ex2_a}
}
\subfigure[]{
\scalebox{0.65}[0.65]{\includegraphics{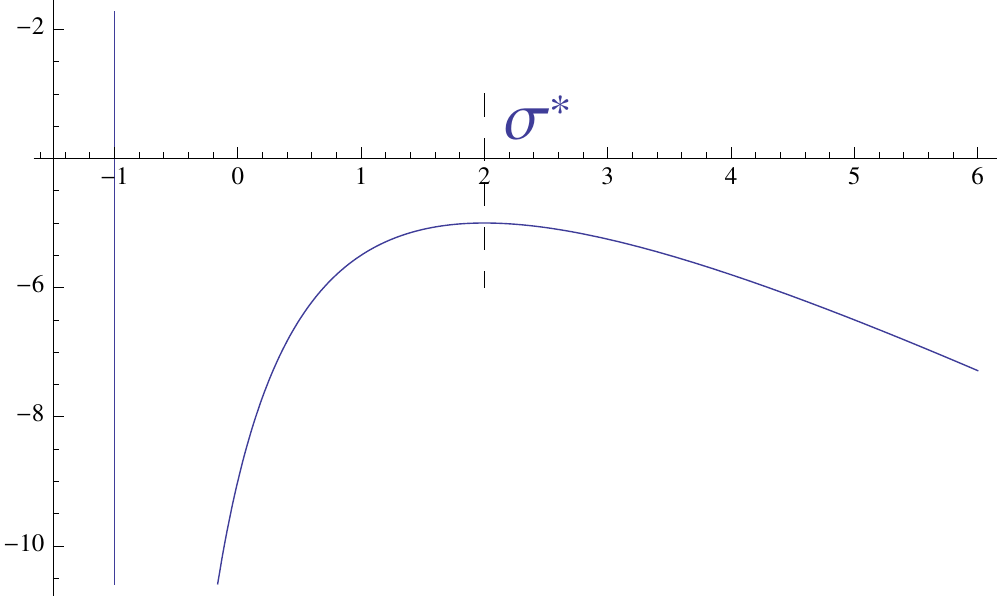}}\label{fig:ex2_b}
}
\caption{ (a): Contours   of the primal function and   boundary of the sphere; (b):  graph of the dual function.}%
\label{fig:ex2}
\end{center}
\end{figure}

{\bf Example 2} The matrix $\QQ$ and radius $r$ are same with that in Example 1 and $\ff$ is changed to
$$
\ff=\bpmat 0\\-3 \epmat,
$$
which is in the same direction of that in Example 1 but has a larger length. We notice that though $\sum_{i=1}^k\hatf_i^2\neq0$ is violated, the condition $\sum_{i=k+1}^n\frac{\hatf_i^2}{(\l_i-\l_1)^2}>r^2$ holds. Thus, it is not in the hard case. There is a critical point in the interior of $\calS_a^+$, which is shown in Figure \ref{fig:ex2_b},  and it is corresponding to the unique global solution of the primal problem, which is the green point in Figure \ref{fig:ex2_a}.

\subsection{Large-size examples}
A hundred of examples are randomly generated, containing fifty examples of the general case and fifty examples of the hard case.
 Both  cases have ten examples
for dimensions of 500, 1000, 2000, 3000 and 5000.
All elements of the coefficients, $\QQ$, $\ff$ and $r$, are integer numbers in $[-100,100]$.
For each example of the hard case, a matrix $\QQ$ with the multiplicity of $\l_1
 =1 $
  is chosen. The corresponding vector $\ff$ is constructed such that $\ff$ is perpendicular to the eigenvector $\UU_1$.
  Then a proper radius $r$ is calculated such that the existence conditions are violated.

Two approaches are used to calculate the value of $\psi(\s)$. One is using decomposition methods to calculate $\GG_a(\s)^{-1}\pp$,
 for which we use the `left division' of Matlab. Another is solving the problem  (\ref{p:xsgivens}),
  for which we use the function `quadprog' of Matlab. For the function `quadprog',
   the tolerance parameter `TolFun' is set as 1e-12. The Matlab is of version 7.13 and runned in the platform with Linux 64-bit system and quad CPUs.

 A perturbation item $\a \UU_1$ is added into the target function for the hard case, and two values of $\a$, 1e-3 and 1e-4, are tried. 
In the main part of the algorithm, the termination tolerance on the value of $\psi(\s)$ is set to be 1e-8.

Results are shown in Table \ref{tb:gene-3}, \ref{tb:gene-4}, \ref{tb:hard-3} and \ref{tb:hard-4}, and they contain the number
of examples which are successfully solved (Succ.Solv.), the distance of the optimal solution to the boundary of the sphere (Dist.Boun.),
 the number of iteration of the Algorithm: Main (Numb.Iter.) and the running time of the algorithm (Runn.Time).
 The values in the columns of Dist.Boun., Numb.Iter. and Runn.Time are averages of the examples successfully solved.
 We compare the results of the algorithm adopting  `left division' and that of the algorithm adopting `quadprog'
 in the same table, where LD denotes left division and QP denotes quadprog.

\renewcommand{\arraystretch}{1}
\begin{table}
\centering
\tbl{General case and $\a=1e-3$.}
{\begin{tabular}{lrrrrrrrrrrr}
\toprule
& \multicolumn{2}{c}{Succ.Solv.}& & \multicolumn{2}{c}{Dist.Boun.} & & \multicolumn{2}{c}{Numb.Iter.} & & \multicolumn{2}{c}{Runn.Time.}\\
\cmidrule(l){2-3} \cmidrule{5-6} \cmidrule{8-9} \cmidrule{11-12}
Dim & LD & QP& & LD & QP & & LD & QP & & LD & QP\\
\hline
500 & 10 & 10 & & 4.716e-09 & 5.245e-09 & & 28.9 & 28.6 & & 0.53 & 1.29\\
1000 & 10 & 10 & & 4.261e-09 & 3.974e-09 & & 27.1 & 27.5 & & 1.67 & 6.25\\
2000 & 10 & 10 & & 3.211e-09 & 3.822e-09 & & 28.2 & 27.8 & & 6.52 & 15.23\\
3000 & 10 & 10 & & 5.674e-09 & 5.221e-09 & & 26.1 & 26.4 & & 20.90 & 72.43\\
5000 & 10 & 10 & & 5.422e-09 & 3.873e-09 & & 28.6 & 28.5 & & 71.68 & 170.34\\
\bottomrule
\end{tabular}}
\label{tb:gene-3}
\end{table}

\renewcommand{\arraystretch}{1}
\begin{table}
\centering
\tbl{General case and $\a=1e-4$.}
{\begin{tabular}{lrrrrrrrrrrr}
\toprule
& \multicolumn{2}{c}{Succ.Solv.}& & \multicolumn{2}{c}{Dist.Boun.} & & \multicolumn{2}{c}{Numb.Iter.} & & \multicolumn{2}{c}{Runn.Time.}\\
\cmidrule(l){2-3} \cmidrule{5-6} \cmidrule{8-9} \cmidrule{11-12}
Dim & LD & QP& & LD & QP & & LD & QP & & LD & QP\\
\hline
500 & 10 & 10 & & 4.532e-09 & 4.464e-09 & & 28.9 & 28.9 & & 0.43 & 1.16\\
1000 & 10 & 10 & & 3.849e-09 & 5.931e-09 & & 27.4 & 27.1 & & 1.47 & 6.08\\
2000 & 10 & 10 & & 2.648e-09 & 2.872e-09 & & 27.9 & 28.5 & & 6.26 & 15.82\\
3000 & 10 & 10 & & 5.299e-09 & 5.137e-09 & & 26.2 & 26.2 & & 20.15 & 73.60\\
5000 & 10 & 10 & & 3.188e-09 & 4.005e-09 & & 28.7 & 28.5 & & 65.71 & 171.92\\
\bottomrule
\end{tabular}}
\label{tb:gene-4}
\end{table}

\renewcommand{\arraystretch}{1}
\begin{table}
\centering
\tbl{Hard case and $\a=1e-3$.}
{\begin{tabular}{lrrrrrrrrrrr}
\toprule
& \multicolumn{2}{c}{Succ.Solv.}& & \multicolumn{2}{c}{Dist.Boun.} & & \multicolumn{2}{c}{Numb.Iter.} & & \multicolumn{2}{c}{Runn.Time.}\\
\cmidrule(l){2-3} \cmidrule{5-6} \cmidrule{8-9} \cmidrule{11-12}
Dim & LD & QP& & LD & QP & & LD & QP & & LD & QP\\
\hline
500 & 10 & 10 & & 4.340e-09 & 6.297e-09 & & 36.0 & 34.9 & & 0.48 & 1.11\\
1000 & 10 & 10 & & 4.253e-09 & 4.904e-09 & & 34.6 & 34.9 & & 1.54 & 3.54\\
2000 & 10 & 10 & & 2.808e-09 & 4.255e-09 & & 35.9 & 35.8 & & 7.15 & 15.11\\
3000 & 9 & 10 & & 5.479e-09 & 4.466e-09 & & 34.0 & 35.0 & & 19.41 & 36.01\\
5000 & 10 & 10 & & 3.755e-09 & 4.705e-09 & & 35.2 & 35.5 & & 74.79 & 121.41\\
\bottomrule
\end{tabular}}
\label{tb:hard-3}
\end{table}

\renewcommand{\arraystretch}{1}
\begin{table}
\centering
\tbl{Hard case and $\a=1e-4$.}
{\begin{tabular}{lrrrrrrrrrrr}
\toprule
& \multicolumn{2}{c}{Succ.Solv.}& & \multicolumn{2}{c}{Dist.Boun.} & & \multicolumn{2}{c}{Numb.Iter.} & & \multicolumn{2}{c}{Runn.Time.}\\
\cmidrule(l){2-3} \cmidrule{5-6} \cmidrule{8-9} \cmidrule{11-12}
Dim & LD & QP& & LD & QP & & LD & QP & & LD & QP\\
\hline
500 & 7 & 9 & & 2.503e-09 & 4.488e-09 & & 39.6 & 40.6 & & 0.51 & 1.36\\
1000 & 9 & 9 & & 3.148e-09 & 4.482e-09 & & 37.4 & 38.3 & & 1.56 & 3.81\\
2000 & 5 & 9 & & 8.668e-09 & 5.785e-09 & & 38.6 & 42.6 & & 7.36 & 17.95\\
3000 & 5 & 10 & & 6.003e-09 & 3.997e-09 & & 38.4 & 40.6 & & 20.43 & 41.06\\
5000 & 8 & 10 & & 4.748e-09 & 2.814e-09 & & 37.8 & 38.8 & & 72.72 & 131.51\\
\bottomrule
\end{tabular}}
\label{tb:hard-4}
\end{table}

We can see that the examples are solved very accurately with error allowance 
 being less than 1e-09, except few instances which are not solved successfully.
  For general cases, all the examples can be solved within no more than 30 iterations, whiles for hard cases,
  the number of iterations is  around 40.
  From the running time, we notice that our method is capable to handle large-size problems in reasonable time.
   The algorithms using 'left division' and 'quadprog' have similar performances in the accuracy and the number of iterations.
    While the one using 'left division' needs much less time than that of the one using 'quadprog'.
    However, the one using 'quadprog' is able to solve more examples successfully.

\section{Conclusion Remarks}\label{se:disc}
We have presented a detailed study on the  quadratic minimization problem with a sphere constraint.
By the canonical duality, this nonconvex optimization is equivalent to a concave maximization dual problem over
a convex domain $\calS^+_a$, which is true also for many other global optimization problems
(see \cite{Gao03perfect,Gao05inequality,Gao06polynomial,Gao07boxnonconvex,gao10fixedcost,Gao10sensor,gao11qp}).
Therefore, the so-called hidden convexity discovered by Ben-Tal and  Teboulle in \cite{BenTal96hidden} is indeed  a special case of the canonical duality theory.
Based on this canonical dual problem,   sufficient and necessary conditions
are  obtained for both general and hard cases.
In order to solve hard case problems,
a perturbation method and the associated algorithm are proposed.
 Numerical results for  large-size examples  demonstrate the efficiency of the proposed approach.
 Combining with the trust region method, the results presented in this paper can be used
   for efficiently  solving general global optimizations.\\

\noindent {\bf Acknowledgements}\\
This  research is supported by US Air Force
Office of Scientific Research under the grant AFOSR FA9550-10-1-0487,  as well as by a  
 grant from the Australian Government under the Collaborative Research Networks
(CRN) program.
The main results of this paper have been announced at the 3rd World Congress of Global Optimization, 
July 9-11, 2013, the Yellow Mountains, China. 

\bibliography{qpsphere}
\bibliographystyle{gOMS}
\end{document}